\documentclass{amsart}

\usepackage{latexsym,amsxtra,amscd,ifthen}

\usepackage{amsfonts}

\usepackage{verbatim}

\usepackage{amsmath}

\usepackage{amsthm}

\usepackage{amssymb}

\numberwithin{equation}{section}

\newcommand{\mattwo}[1]{\left[\begin{array}{rr} #1 \end{array}\right]}

\newcommand{\matthree}[1]{\left[\begin{array}{rrr} #1 \end{array}\right]}

\theoremstyle{plain}

\newtheorem{theorem}{Theorem}[section]

\newtheorem{lemma}[theorem]{Lemma}

\newtheorem{proposition}[theorem]{Proposition}

\newtheorem{conjecture}[theorem]{Conjecture}

\theoremstyle{definition}

\newtheorem{definition}[theorem]{Definition}

\newtheorem{example}[theorem]{Example}

\newtheorem{hypotheses}[theorem] {Hypotheses}
\newtheorem{question}[theorem]{Question}

\makeatletter              

\let\c@equation\c@theorem  

\makeatother

\newcommand{\rank}{\operatorname{rank}}

\DeclareMathOperator{\hdet}{hdet}

\DeclareMathOperator{\Ext}{Ext}

\DeclareMathOperator{\Aut}{Aut}

\DeclareMathOperator{\GKdim}{GKdim}

 \DeclareMathOperator{\Hom}{Hom}

\newcommand{\mf}{\mathfrak}

\newcommand{\mc}{\mathcal}

\newcommand{\be}{\begin{enumerate}}

\newcommand{\ee}{\end{enumerate}}

\newcommand{\bq}{\begin{eqnarray*}}

\newcommand{\eq}{\end{eqnarray*}}

\newcommand{\bqn}{\begin{eqnarray}}

\newcommand{\eqn}{\end{eqnarray}}

\begin{document}

\title[Invariant Theory]
{invariant theory of artin-schelter regular algebras: a survey}

\author{Ellen E. Kirkman}

\address{ Department of Mathematics,
P. O. Box 7388, Wake Forest University, Winston-Salem, NC 27109}
\email{kirkman@wfu.edu}

\begin{abstract}
This is survey of results that extend notions of the classical invariant theory of linear actions by finite groups on $k[x_1, \dots, x_n]$ to the setting of finite group or Hopf algebra $H$ actions on an Artin-Schelter regular algebra $A$.  We investigate when $A^H$ is AS regular, or AS Gorenstein, or a ``complete intersection" in a sense that is defined.  Directions of related research are explored briefly.
\end{abstract}


\maketitle

\bigskip

\setcounter{section}{-1}

\section{Introduction}
\label{xxsec0}
The study of invariants of finite groups acting on a commutative polynomial ring $k[x_1, \dots, x_n]$ has played a major role in the development of commutative algebra, algebraic geometry, and representation theory.  This paper is a survey of more recent work  that extends these techniques to a noncommutative setting.   We will be particularly concerned with algebraic properties of the subring of invariants.

We begin with some basic definitions.  Throughout we let $k$ be an algebraically closed field of characteristic zero and $A$ be a $k$-algebra.   A $k$-algebra $A$ is said to be {\it connected graded} if
$A = k \oplus A_1 \oplus A_2 \oplus \cdots$
with $A_i \cdot A_j \subseteq A_{i+j}$ for
all $i,j \in \mathbb{N}$; we denote the trivial module of a connected graded algebra by $k$. Throughout we assume that a connected graded algebra $A$ is generated in degree 1. The {\it Hilbert series} of $A$ is defined
to be the formal power series $H_A(t) = \sum_{i \in \mathbb{N}} (\dim_k A_i) t^i$.  The Gelfand-Kirillov dimension of an algebra $A$ is
denoted by $\GKdim A$;  it is related to the rate of growth
in the dimensions of the graded pieces $A_n$ of $A$ (see \cite{KL}).
The commutative polynomial ring $k[x_1, \dots, x_n]$
has Gelfand-Kirillov dimension $n$. 
In our noncommutative setting we replace the commutative polynomial ring $k[x_1, \dots, x_n]$ with an Artin-Schelter regular algebra, defined as follows.

\begin{definition}
\label{zzdef1.1}
Let $A$ be a connected graded algebra.
We call $A$ {\it Artin-Schelter Gorenstein} (or {\it AS Gorenstein}
for short) {\it of dimension $d$} if the following conditions hold:
\begin{enumerate}
\item[(a)]
$A$ has injective dimension $d<\infty$ on the left and
on the right,
\item[(b)]
$\Ext^i_A(_Ak,_AA)=\Ext^i_{A}(k_A,A_A)=0$ for all $i\neq d$, and
\item[(c)]
$\Ext^d_A(_Ak,_AA)\cong \Ext^d_{A}(k_A,A_A)\cong k(-l)$ for some $l$ (where
$l$, the shift in the grading, is called the {\it AS index} of $A$).
\end{enumerate}
If, in addition,
\begin{enumerate}
\item[(d)]
$A$ has finite global dimension, and
\item[(e)]
$A$ has finite Gelfand-Kirillov dimension,
\end{enumerate}
then $A$ is called {\it Artin-Schelter regular}
(or {\it AS regular} for short) {\it of dimension $d$}.
\end{definition}

Note that polynomial rings $k[x_1,\dots, x_n]$
for $n\geq 1$, with $\deg x_i=1$, are AS
regular of dimension $n$, and they are the only
commutative AS regular algebras.  Hence
AS regular algebras are natural
generalizations of commutative polynomial rings.

In some cases we are able to prove stronger results for a special class of
AS regular algebras that we call quantum polynomial rings.
\begin{definition}
\label{quantumpolynomial} Let $A$ be a connected graded
algebra. If $A$ is a noetherian, AS  regular graded
domain of global dimension $n$ and $H_A(t)=(1-t)^{-n}$,
then we call $A$ {\it a quantum polynomial ring
of dimension $n$}.
\end{definition}

By \cite[Theorem 5.11]{Sm2}, a quantum polynomial ring is
Koszul and hence it is generated in degree 1. The
GK-dimension of a quantum polynomial ring of global
dimension $n$ is $n$. 

Artin-Schelter regular algebras of dimension 3 were classified in
\cite{ASc, ATV}. They occur in two families: (1) the quantum polynomial
algebras having 3 generators and 3 quadratic relations (that include the {\it 3-dimensional Sklyanin algebra}), and (2) algebras having 2 generators and 2 cubic relations (that include the noetherian graded {\it down-up algebras}, which will be discussed in Section 3).  There are many examples of AS regular algebras of higher dimensions,  but their
classification has not been completed (current research centers on dimension 4, which appears to be much more complex).  The invariant theory of AS regular algebras may become richer as higher dimensional AS regular algebras are discovered and classified.

We consider finite groups $G$ of graded automorphisms acting on $A$, and throughout we denote the graded automorphism group of $A$ by $\Aut(A)$.  We note that any $n \times n$ matrix acts naturally on a commutative polynomial ring in $n$ indeterminates.  However, in order for a linear map defined on the degree 1 piece of  a noncommutative algebra $A$ to give a well-defined graded homomorphism of $A$, one must check that the map preserves the ideal of relations; for example, the transposition of $x$ and $y$ is an automorphism of $k_q[x,y]$, the skew polynomial ring with the relation $yx=qxy$, if and only if $q = \pm 1$.  

Typically there is a limited supply of graded automorphisms of $A$, so we introduce further noncommutativity by allowing actions on $A$  by a finite dimensional Hopf algebra $H$. Throughout we adopt the usual notation for the Hopf structure of a Hopf algebra
$H$, namely $(H, m, \epsilon, \Delta, u, S)$ (see \cite{Mon}, our basic reference for Hopf algebra actions), and we denote the $H$-action on $A$ by $\cdot : H
\otimes A \rightarrow A$.  All groups $G$ we consider are finite, and all Hopf algebras $H$ are finite dimensional.

\begin{definition}\
\label{def1.3}
\begin{enumerate}
\item Let $H$ be a Hopf algebra $H$ and $A$ be a $k$-algebra. We say
that {\it $H$ acts on $A$} (from the left), or $A$ is a {\it left
$H$-module algebra}, if $A$ is a left $H$-module, and  for all $h \in H$,
$h \cdot (ab) = \sum (h_1 \cdot a)(h_2 \cdot b)$  for all $a,b \in A$, where $\Delta(h) =
\sum h_1 \otimes h_2$ (using the usual Sweedler notation convention),
 and $h \cdot 1_A = \epsilon(h)
1_A$. 
\item The {\it invariant subring} of
such an action is defined to be
$$A^H = \{ a \in A ~|~ h \cdot a = \epsilon(h) a, ~\forall ~h \in H \}.$$
\end{enumerate}
\end{definition}
When the Hopf algebra $H=k[G]$ is the group algebra of a finite group $G$, the usual coproduct on $k[G]$ is $\Delta(g) = g \otimes g$, and the Hopf algebra action on $A$ means that $g$ acts as a homomorphism on elements of $A$.  Similarly, since $\epsilon(g) = 1$, the invariant subring $A^{k[G]} = A^G$ is the usual subring of invariants.

Sometimes it is more convenient to view a left $H$-module algebra $A$ as a right $H^\circ$-comodule algebra, where $H^\circ$ is the Hopf dual of $H$.  When $H$ is finite dimensional the Hopf dual is the vector space dual, and a left $H$-module can be viewed as a right $H^\circ$-module.  To be a right $K$-comodule algebra over a Hopf algebra $K$ we require a coaction $\rho: A \rightarrow A \otimes K$ with properties: $\rho(1_A) = 1_A \otimes 1_K$ and
$\rho(ab) = \rho(a) \rho(b)$ for all $a,b \in A$.  The coinvariant subring of such a coaction is defined to be
$$A^{co K} = \{ a \in A | \rho(a) = a \otimes 1_K \}.$$
It follows (\cite[Lemma 6.1 (c)]{KKZ2}) that $A^H = A^{co H^\circ}$.

Throughout we require that the Hopf algebra $H$ acts on $A$ so that 
\begin{hypotheses} \label{standing}\
\begin{itemize}
\item $A$ is an $H$-module algebra, 
\item the grading on $A$ is preserved, and 
\item the action of $H$ on $A$ is inner faithful.
\end{itemize}
\end{hypotheses} 
 The inner faithful assumption guarantees that the Hopf algebra action is not actually an action over a homomorphic image of $H$ that might be a group algebra (in which case the action is actually a group action).
\begin{definition}\cite{BB}.
\label{def1.4} Let $M$ be a left $H$-module. We say that $M$ is an {\it
inner faithful} $H$-module, or $H$ \emph{acts inner faithfully} on $M$,
if $IM\neq 0$ for every nonzero Hopf ideal $I$ of $H$.
Dually, let $N$ be a right $K$-comodule.  We say that $N$ is an {\it inner faithful} $K$-comodule, or that $K$ {\it acts inner faithfully} on $N$, if for any proper Hopf subalgebra $K' \subsetneq K$, $\rho(N)$ is not contained in $N \otimes K'$.
\end{definition}

\begin{lemma} \cite[Lemma 1.6(a)]{CKWZ1}. Let $H$ be a finite dimensional Hopf algebra, $K = H^\circ$, and $U$ be a left $H$-module.  Then $U$ is a right $K$-comodule, and the $H$-action on $U$ is inner faithful if and only if the induced $K$-coaction on $U$ is inner faithful.
\end{lemma}

Despite the hope that Hopf algebra actions would provide many new actions on $A$, we note the result of Etingof and Walton, that states that if $H$ is a semisimple Hopf algebra acting inner faithfully on a commutative domain, then $H$ is a group algebra \cite{EW1};  this result does not require our assumptions that the action of $H$ on $A$ preserves the grading on $A$.  Other results showing that  Hopf actions must factor through group actions are mentioned in Section 4. 
There are, however, interesting actions by pointed Hopf algebras on domains and fields \cite{EW2}.   Moreover, there are subrings of invariants that occur under Hopf actions that are not invariants under group actions, warranting continued study of these actions.

In the first section we discuss the question of when $A^H$ is AS regular, extending the Shephard-Todd-Chevalley Theorem to our noncommutative context.  In the second section we discuss the question of when $A^H$ is AS Gorenstein; these results include extending Watanabe's Theorem that $k[x_1, \dots, x_n]^G$ is Gorenstein when $G$ is a finite subgroup of ${\rm SL}_n(k)$, and Felix Klein's classification of the invariants of finite subgroups of ${\rm SL}_2(k)$ acting on $k[x,y]$.  In the third section we discuss the question of when $A^G$ is a ``complete intersection" (and what a complete intersection might mean in this context); here our goal is to extend results of Gordeev, Kac, Nakajima, and Watanabe to our noncommutative setting.  In the final section we briefly discuss some related directions of research.   
\section{Artin-Schelter regular subrings of invariants}

\label{xxsec1}
A result attributed to Carl F. Gauss states that the subring of invariants  of the commutative polynomial ring $k[x_1, \dots, x_n]$  under the action of the symmetric group $S_n$, permuting  the indeterminates, is generated by the $n$ elementary symmetric polynomials.  Further, the symmetric polynomials are algebraically independent, so that the subring of invariants is also a polynomial ring.  This result raised the question:  for which finite groups $G$ acting on $k[x_1, \dots, x_n]$ is the fixed subring a polynomial ring?  This question was answered in 1954 for $k$ an algebraically closed field of characteristic zero by G. C.  Shephard and J. A. Todd \cite{ShT}, who classified the complex reflection  groups and produced their invariants.  Shortly afterward, C. Chevalley \cite{C} gave a more abstract argument that showed that for real reflection groups $G$, the fixed subring $k[x_1, \dots, x_n]^G$ is a polynomial ring, and J.-P. Serre \cite{S} showed that Chevalley's argument could be used to prove the result for all unitary reflection groups.  A. Borel's history \cite[Chapter VII]{B} provides details on the origins of the ``Shephard-Todd-Chevalley Theorem".  Invariants in a commutative polynomial ring under the action of a reflection group do \underline{not} always form a polynomial ring in characteristic p (see \cite[Example 3.7.7]{DK} (and following) for a discussion, including an example of Nakajima \cite{N1}).  However, in characteristic p, to obtain a polynomial subring of invariants it is necessary (but not sufficient) that the group be a reflection group (see \cite[proof of Theorem 7.2.1]{Be} for a proof that works in any characteristic).   In characteristic zero the necessity of $G$ being a reflection group follows from the sufficiency by considering the normal subgroup of reflections in the group.
\begin{theorem}[{\bf Shephard-Todd-Chevalley Theorem}]  \cite{ShT} {\rm and} \cite{C}. For $k$ a field of characteristic zero, the subring of invariants $k[x_1, \dots, x_n]^G$ under a finite group $G$ is a polynomial ring if and only if $G$ is generated by reflections.  
\end{theorem}
In this context a linear map $g$ on a vector space $V$, where $g$ is of finite order and hence is diagonalizable, is called a {\it reflection of $V$} if all but one of the eigenvalues of $g$ are 1 (i.e. dim $V^g = $ dim $V - 1$) (we note that sometimes the term ``reflection" is reserved for the case that $g$ has real eigenvalues (so that the single non-identity eigenvalue is -1) and the term ``pseudo-reflection" is used for the case that the single non-identity eigenvalue is a complex root of unity).  We will call a graded automorphism that is a reflection (in this sense) of $V=A_1$, the $k$-space of elements of degree 1,  a ``classical reflection".

In the setting of the Shephard-Todd-Chevalley Theorem the fixed subring $A^G$ is isomorphic to the original ring $A$ (both being commutative polynomial rings), and early noncommutative generalizations of the Shephard-Todd-Chevalley Theorem focused on this  property.  S.P. Smith \cite{Sm1} showed that if $G$ is a finite group acting on the first Weyl algebra
$A= A_1(k)$ then $A^G$ is isomorphic to $A$ if and only if $G = \{ 1 \}$. Alev and Polo \cite{AP} extended Smith's result to the higher Weyl algebras, and  showed further that if ${\mathfrak g}$ and ${\mathfrak g}'$ are two
semisimple Lie algebras, and $G$ is a finite group of
algebra automorphisms of the universal enveloping algebra $U({\mathfrak g})$ such that
$U({\mathfrak g})^G \cong U({\mathfrak g}')$, then $G$
is trivial and ${\mathfrak g}\cong {\mathfrak g}'$.  The preceding results were attributed to the ``rigidity" of noncommutative algebras, and these early results suggested that there is no noncommutative analogue of the Shephard-Todd-Chevalley Theorem.  However, as we shall see,  in the case that the algebra $A$ is graded there are other ways to generalize the Shephard-Todd-Chevalley Theorem.

We begin with an illustrative noncommutative example.
\begin{example}
\label{classicalref}
Let $A=k_{-1}[x,y]$ be the skew polynomial ring with the relation $yx=-xy$; the ring $A$ is an AS regular algebra of dimension 2.  Let $G= \langle g \rangle $ be the cyclic group generated by the graded automorphism $g$, where $g(x) = \lambda_n x$ and $g(y) = y$, for
 $\lambda_n$  a primitive n-th root of unity.  The linear map $g$ acting on $V=A_1$ is a classical reflection.  The fixed ring  ${A}^{G} = k \langle x^n, y \rangle = k_{(-1)^n}[x,y]$ is isomorphic to $A$ when $n$ is odd, but not when $n$ is even (when it is a commutative polynomial ring).  However, $A^G$ is
AS regular for all $n$.  This example suggests that a reasonable generalization of the Shephard-Todd-Chevalley Theorem is that $G$ should be thought of as a ``reflection group" when $A^G$ is AS regular, rather than when $A^G$ is isomorphic to $A$.  When $A$ is a commutative polynomial ring these two conditions are the same condition:  when $A$ is a commutative polynomial ring, $A^G$ is AS regular if and only if $A^G \cong A$.
\end{example}

\begin{definition} We call a finite group $G$  of graded automorphisms of an AS regular algebra $A$ a {\it reflection group for $A$} if the fixed subring $A^G$ is an AS regular algebra.
\end{definition}

 In our terminology the classical reflection groups are reflection groups for\linebreak  $k[x_1, \dots, x_n]$.  The next example demonstrates that a classical reflection group will not always produce an AS regular invariant subring when acting on some other AS regular algebra.

\begin{example}
\label{transposition}  Again, let $A=k_{-1}[x,y]$ be the skew polynomial ring with the relation $yx=-xy$.
The transposition $g$ that interchanges $x$ and $y$ induces a graded automorphism of $A$, and $g$ generates the symmetric group $S_2$, which is a classical reflection group.  One set of generators for the fixed ring $A^{\langle g \rangle}$ is $x+y$ and $x^3 +y^3$.  Here $xy$ is not fixed, as it is in the commutative case, and the invariant $x^2 + y^2 = (x + y)^2$ is not a generator.  The generators $x+y$ and $x^3 +y^3$ are not algebraically independent, and the algebra they generate is not an AS regular algebra.  However, as we shall see (Example \ref{transagain}), the fixed ring is AS Gorenstein, and it can be viewed as a hypersurface in an AS regular algebra of dimension 3.
\end{example}

Example \ref{transposition} shows that when $A$ is noncommutative  we need a different notion of  ``reflection" to obtain an AS regular fixed subring, as the eigenvalues of the linear map $g$ no longer control the AS regularity of the fixed subring.  Our results suggest that it is the trace function, defined below, that determines whether a linear map is a ``reflection".  In Example \ref{toshowtraces} below we will show that  the classical reflection $g$ in Example \ref{transposition} is not a reflection in this new sense, while the automorphism $g$ of Example \ref{classicalref} is a reflection (under both the classical definition and our new definition).

\begin{definition} Let $A$ be a graded algebra with $A_j$ denoting the elements of degree $j$.  The {\it trace function} of a graded automorphism $g$ acting on $A$ is defined to be the formal power series
$$Tr_A(g,t) = \sum_{j=0}^{\infty} trace(g|_{A_j})\; t^j,$$
where  $trace(g|_{A_j})$ is the usual trace of the linear map $g$ restricted to $A_j$.
\end{definition}
In this setting, by \cite[Theorem 2.3(4)]{JiZ} $Tr_A(g,t)$ is a rational function of the form
$1/e_g(t)$, where $e_g(t)$ is a polynomial in $k[t]$, and the zeroes of $e_g(t)$ are all roots of unity \cite[Lemma1.6(e)]{KKZ1} (in the case $A=k[x_1, \dots, x_n]$, the roots of $e_g(t)$ are the inverses of the eigenvalues of $g$).  The next proposition shows that trace functions can be used in computing fixed subrings, giving a version of the classical Molien's Theorem.  Knowing the Hilbert series of the subring of invariants is very useful in computing the subring of invariants. (In the Hopf algebra action case of the theorem below, see \cite[Definition 2.1.1]{Mon} for the definition of the integral).
\begin{proposition}[{\bf{Molien's Theorem}}] \label{Molien}  \cite[Lemma 5.2]{JiZ}, \cite[Lemma 7.3]{KKZ2}.  The Hilbert series of the fixed subring $A^G$ is 
$$H_{A^G}(t) = \frac{1}{|G|} \sum_{g \in G} Tr_A(g,t).$$ Similarly, for a semisimple Hopf algebra $H$ acting on $A$ with integral $\int$ that has $\epsilon(\int) = 1 \in k$,  the Hilbert series of the fixed subring $A^H$ is $H_{A^H}(t) = Tr_A(\int,t)$.
\end{proposition}

In our setting it is the order of the pole of $Tr_A(g,t)$ at $1$, rather than the eigenvalues of $g$, that determines whether $g$ is a ``reflection".
\begin{definition} \cite[Definition 1.4]{KKZ1}.  Let $A$ be an AS regular algebra of $\GKdim A$ $=  n$.
We call a graded automorphism $g$ of $A$ a {\it reflection of $A$} if the trace function of $g$ has the form
$$Tr_A(g,t) = \frac{1}{(1-t)^{n-1} q(t)} \mbox{ where } q(1) \neq 0.$$
\end{definition} 
We note that in \cite{KKZ1} we used the term ``quasi-reflection" to distinguish our use of reflection from the usual notion of reflection.  Here we will use the term ``reflection" as defined above, and refer to ``classical reflection" when we are referring to a reflection defined in terms of its eigenvalues.
The following examples can be used to justify our definition of reflection, as we want a group generated by ``reflections" to have a fixed subring that is AS regular (i.e. to be a reflection group for $A$).
\begin{example}
\label{toshowtraces}
Let $A= k_{-1}[x,y]$, an AS regular algebra of dimension 2, and in each case let $G= \langle g \rangle $ be the cyclic group generated by the graded automorphism $g$, expressed as a matrix acting on the vector space $A_1 = kx \oplus ky$.  Let $\lambda_n$ be a primitive n-th root of unity.
\begin{enumerate}

\item  As in Example \ref{classicalref}, let ${\displaystyle  {g = \mattwo{ \lambda_n & 0\\ 0 & 1}}}$;  the automorphism $g$ is a classical reflection of $A$.   The trace function is ${\displaystyle Tr_{A}(g ,t) = \frac{1} {(1-t)(1- \lambda_n t)}}$, so that $g$ is a reflection of $A$ under our definition.  Furthermore   ${A}^{G} = k \langle x^n, y \rangle = k_{(-1)^n}[x,y]$ is AS regular for all $n$.
 
\item As in Example \ref{transposition}, let ${\displaystyle  {g=  \mattwo{0 & 1\\1 & 0}}}$.  The trace function is  ${\displaystyle Tr_{{A}}({g},t) = \frac{1}{1+t^2}}$.  As we noted  earlier, ${A}^G$ is not AS regular. The automorphism $g$ is a classical reflection, but it is not a reflection of $A$ in our sense (the order of the pole at 1 is not 2-1=1, and the fixed subring is not AS  regular).

\item  Let ${\displaystyle {g = \mattwo{0 & -1\\1 & 0}}}$; the automorphism $g$ is not a classical reflection.  The trace function is  ${\displaystyle Tr_{{A}}({g},t) = \frac{1}{{(1-t)}(1+t)}}$ and  ${A}^{G} = k[x^2 + y^2, xy]$ is a commutative polynomial ring so is  AS regular. The automorphism $g$ is a reflection of $A$ in our sense, and we called it a ``mystic reflection" to distinguish it from a classical reflection (a reflection such as in (1)).\\
\end{enumerate}
\end{example}
Proposition \ref{Molien} is used in computing  Example \ref{toshowtraces}.  For example, in (3) Molien's Theorem shows us that
$$H_{{A}^{G}}(t) = \frac{1}{4(1-t)^2}  + \frac{2}{4(1-t^2)}  + \frac{1}{4(1+t)^2}  = \frac{1}{(1-t^2)^2},$$
hence the two invariants we have found generate the invariant subring, because they are fixed and the ring they generate has this Hilbert series.

We have shown that for $A$ a quantum polynomial ring, there are only two kinds of reflections of $A$:  classical reflections and new reflections (as in Example \ref{toshowtraces} (3)) that we call ``mystic reflections".
\begin{theorem}\cite[Theorem 3.1]{KKZ1}.
\label{xxthm3.1} Let $A$ be a quantum polynomial ring of
global dimension $n$. If $g\in \Aut(A)$ is a reflection of $A$
of finite order, then $g$ is in one of the following two
cases:
\begin{enumerate}
\item
There is a basis of $A_1$, say $\{b_1,\cdots,b_n\}$,
such that $g(b_j)=b_j$ for all $j\geq 2$ and $g(b_1)
=\lambda_n b_1$ for $\lambda_n$ a root of unity.  Hence $g|_{A_1}$ is a classical reflection.
\item
The order of $g$ is $4$ and there is a basis of $A_1$,
say $\{b_1,\cdots,b_n\}$,
such that $g(b_j)=b_j$ for all $j\geq 3$ and $g(b_1)
=i\; b_1$ and $g(b_2)=-i\; b_2$ (where $i^2=-1$).  We call such a reflection a {\it mystic reflection}.
\end{enumerate}
\end{theorem}

It is quite possible that other AS regular algebras have different kinds of reflections, as AS regular algebras have been completely classified only in dimensions $\leq 3$.

We conjecture the following generalization of the Shephard-Todd-Chevalley Theorem.
\begin{conjecture}\label{STCconjecture}
Let $A$ be an AS regular algebra and $G$ be a finite group of graded automorphisms of $A$.  Then $A^G$ is AS regular if and only if $G$ is generated by reflections of $A$.
\end{conjecture}
We have proved a number of partial results that support this conjecture.
First we note that, due to the following theorem, only one direction of this conjecture needs to be proved when $A$ is noetherian.
\begin{theorem}\cite[Proposition 2.5(b)]{KKZ3}. Let $A$ be a noetherian AS regular algebra and suppose that for every finite group of graded automorphisms of $A$ that is generated by reflections of $A$ it follows that $A^G$ is AS regular.  Then Conjecture \ref{STCconjecture} is true.
\end{theorem}
Although we have not proved that if $A^G$ is AS regular, then $G$ is generated by reflections, we have shown that $G$ must contain at least one reflection.
\begin{theorem}\cite[Theorem 2.4]{KKZ1}.
\label{xxthm2.4}
Let $A$ be noetherian and AS regular, and let
 $G$ be a finite group of graded automorphisms of $A$. If $A^G$ has
finite global dimension, then $G$ contains a reflection of $A$.
\end{theorem}
We have shown that a number of algebras have no reflections (hence no AS regular fixed algebras).  These algebras include: non-PI Sklyanin algebras \cite[Corollary 6.3]{KKZ1}, homogenizations of the universal enveloping algebra of a finite dimensional Lie algebras $\mathfrak{g}$ (\cite[Lemma 6.5(d)]{KKZ1}), and down-up algebras (or any noetherian AS regular algebra of dimension 3 generated by two elements of degree 1) (\cite[Proposition 6.4]{KKZ1}).  We will say more about down-up algebras in Section 3.

The skew polynomial ring $k_{p_{ij}}[x_1, \dots, x_n]$ is defined to be the $k$-algebra generated by $x_1, \dots, x_n$ with relations $x_jx_i = p_{ij}x_ix_j$ for all $1 \leq i < j \leq n$ and $p_{ii} =1$.
\begin{theorem}\label{skewed} \cite[Theorem 5.5]{KKZ3}. Let $A =k_{p_{ij}}[x_1, \dots, x_n]$, and let $G$ be a finite group of graded automorphisms of $A$.  Then $A^G$ has finite global dimension if and only if $G$ is generated by reflections of $A$ (in which case $A^G$ is again a skew polynomial ring).
\end{theorem}
Theorem \ref{skewed} has been proved using different techniques by Y. Bazlov and A. Berenstein in \cite{BB2}, and we will say more about their results shortly.  
\begin{theorem}\cite[Theorem 6.3]{KKZ3}.
Let $A$ be a quantum polynomial ring and let $G$ be a finite abelian group of graded automorphisms of $A$.  Then $A^G$ has finite global dimension if and only if $G$ is generated by reflections of $A$.
\end{theorem}

In their seminal paper Shephard and Todd  classified the (classical) reflection groups, i.e. the reflection groups for  $k[x_1, \dots, x_n]$.  When $A$ is a noncommutative AS regular algebra, whether or not a group is a reflection group depends upon the algebra $A$ on which the group acts, and groups different from the classsical reflection groups can occur as reflection groups for some noncommutative AS regular algebra.  
  We present two examples of reflection groups on $k_{p_{ij}}[x_1, \dots, x_n]$ with $p_{ij} = -1$ for all $i \neq j$.
\begin{example}\cite[Example 7.1]{KKZ3}.  Let $G$ be the group generated by the mystic reflections
$$ g_1 = \begin{pmatrix}
0 & -1 & 0\\
1 & 0 & 0\\
0 & 0 &1\\
\end{pmatrix}
\text {and }
g_2 = \begin{pmatrix}
1 & 0 & 0\\
0 & 0 & -1\\
0 & 1 & 0\\
\end{pmatrix}.$$
$G$ is the rotation group of the cube, and is isomorphic to the symmetric group $S_4$.  $G$ acts on $k_{-1}[x_1, x_2, x_3]$ where $p_{ij} = -1$ for $i \neq j$, with fixed ring the commutative polynomial ring $k[x_1^2+x_2^2+x_3^2, \; x_1x_2x_3, \; x_1^4+x_2^4+x_3^4]$.  Hence under this representation (but not the permutation representation) $S_4$ is a reflection group for $k_{-1}[x_1, x_2, x_3]$.
\end{example}
\begin{example} \cite[Example 7.2]{KKZ3}. The binary dihedral groups
$G  = BD_{4 \ell}$ are generated by the mystic reflections
$$ g_1 = \begin{pmatrix}
	0& -\lambda^{-1}\\
\lambda & 0
\end{pmatrix} \text{ and }
{g_2 = \begin{pmatrix}
0 & -1\\
1 & 0
\end{pmatrix}}$$
for $\lambda$ a primitive $2 \ell$-th root of unity. These groups act on $A=k_{-1}[x,y]$  with fixed ring $A^G = k[x^{2 \ell}+ y^{2 \ell}, xy]$, a commutative polynomial ring.  Hence the binary dihedral groups are reflection groups for $k_{-1}[x,y]$.
When $\ell= 2$, $G$ is the quaternion group of order 8, which is not a classical reflection group. 
\end{example}

 Other examples of reflections groups for $k_{-1}[x_1, \dots, x_n]$ include the infinite family $M(n,\alpha, \beta)$ (\cite[Section 7]{KKZ3}); there are infinite families of these groups that are not isomorphic as groups to classical reflection groups.  Bazlov and Berenstein found this same class of groups occurring in their work \cite{BB1} related to Cherednik algebras, and in \cite{BB2} gave a different proof of Theorem \ref{skewed} by introducing a non-trivial correspondence between reflection groups $G$ for $k_{p_{ij}}[x_1, \dots, x_n]$ and classical reflection groups $G'$.  In particular, they showed that for $G$ a reflection group for $k_{q_{ij}}[x_1, \dots, x_n]$  the group algebra $k[G]$ is isomorphic to the group algebra $k[G']$ (as algebras) for $G'$ a classical reflection group, even though $G$ and $G'$ are not isomorphic as groups.

Actions of noncocommutative Hopf algebras also can produce AS regular fixed subrings. Hence we  expand our notion of reflection group to include ``reflection Hopf algebras" for a given AS regular algebra.

\begin{definition} We call a Hopf algebra $H$ a {\it reflection Hopf algebra for $A$} if  $A^H$ is AS regular.
\end{definition}

The group algebra of a classical reflection group is a reflection Hopf algebra for $k[x_1, \dots,x_n]$. Next we present an example of a Hopf algebra that is not commutative or cocommutative but has an AS regular ring of invariants.
\begin{example}\cite[Section 7]{KKZ2}. 
The smallest dimensional semisimple Hopf algebra $H$ that is not isomorphic to a group algebra or its dual is the 8-dimensional semisimple algebra $H_8$, defined by Kac and Paljutkin \cite{KP} (see  also  \cite{Ma1}).

As an algebra {$H_8$} is generated by {$x, y, z$} with the following relations:
$$x^2 = y^2 =1, \;\; xy=yx,\;\; zx=yz,$$
$$ zy=xz,\;\; z^2= \frac{1}{2}(1+x+y-xy). $$
The coproduct, counit and antipode are given as follows:
$$\Delta(x) = x\otimes x, \;\;\; \Delta(y)=y\otimes y,$$
$$ \Delta(z) = \frac{1}{2}(1\otimes 1 + 1\otimes x + y\otimes 1 - y\otimes x)(z\otimes z), $$
$$ \epsilon(x) = \epsilon(y) = \epsilon(z) =1, \quad S(x)=x^{-1},\; S(y)=y^{-1},\; S(z)=z.$$
 ${H_8}$ has a unique irreducible 2-dimensional representation on
$k u \oplus k v$ given by
{$$ x \mapsto
\begin{pmatrix}-1 & 0\\ 0 & 1
\end{pmatrix}, \quad
\;
y \mapsto \begin{pmatrix}1 & 0\\ 0& -1
\end{pmatrix}, \quad
\;
z \mapsto \begin{pmatrix}0 & 1 \\1& 0
\end{pmatrix}.$$}

\begin{enumerate}
\item Let ${A = k\langle u, v \rangle / \langle u^2-v^2 \rangle}$ ($A$ is isomorphic to $k_{-1}[u,v]$). Then $H_8$ acts on $A$ and the fixed subring is 
${A}^{H_8} = k[u^2, (uv)^2 - (vu)^2],$
a commutative polynomial ring, and so  ${H_8}$ is a reflection Hopf algebra for ${A}$.
\item  Let ${A =k\langle u, v \rangle / \langle vu \pm iuv \rangle}$.  Then $H_8$ acts on $A$ and the fixed subring is
${A}^{H_8} = k[u^2v^2, u^2 + v^2]$,
a commutative polynomial ring. Hence ${H_8}$ is  a reflection Hopf algebra for ${A}$.
\end{enumerate}
 
\end{example}
Furthermore, actions of non-semisimple Hopf algebras can produce AS regular fixed algebras.
\begin{example}\cite[Section 3.2.1]{All}.
The Sweedler algebra ${H(-1)}$ is generated by ${g}$ and ${x}$ with algebra relations:
$$  g^2=1,\;\; x^2 = 0, \;\; xg = - gx,$$
and coproduct, counit, and antipode:
$$ \Delta(g) = g\otimes g  \;\;\;\; \Delta(x) =   g\otimes x + x\otimes 1, $$
$$  \epsilon(g) = 1, \; \epsilon(x) = 0 \quad S(g)=g,\; S(x)=-g x.$$
 Then ${H(-1)}$ acts on the commutative polynomial algebra ${k[u,v]}$ as
{$$ x \mapsto
\begin{pmatrix} 0 & 1\\ 0 & 0
\end{pmatrix}, \quad
\;
g \mapsto \begin{pmatrix}1 & 0\\ 0& -1
\end{pmatrix}
$$}
with fixed subring ${k[u,v]}^{H(-1)} = k[u,v^2],$ a commutative polynomial ring.  Hence $H(-1)$ is a reflection Hopf algebra for $k[u,v]$.
\end{example}

In work in progress we have shown that (using the notation of \cite{Ma2}) the Hopf algebras $H= A_{4m}$ (for $m$ odd)  and $H= B_{4m}$ are reflection Hopf algebras for $A=k_{-1}[x,y]$; in both cases as an algebra (but not as a Hopf algebra) $H$ is isomorphic to $k[D_{4m}]$, the group algebra of the dihedral group of order $4m$, a classical reflection group.  Further $A_{12}$ acts on a 3-dimensional (non-PI) AS regular algebra $A$ with (non-PI) regular fixed ring.  These examples, along with the examples of group algebras acting on skew-polynomial algebras, suggest that the algebra structure of $H$ and its relation to the group algebra of a classical reflection group may be related to  conditions that guarantee that the Hopf algebra is a reflection Hopf algebra.  We also have some examples of commutative Hopf algebras that are reflection Hopf algebras; in this case the algebra structure of $H$ is not informative.

While we have made Conjecture \ref{STCconjecture} on the properties of a group $G$ that make it a reflection group for $A$, we have made no conjectures on the properties of a general Hopf algebra $H$ that make it a reflection Hopf algebra for $A$.  One can still take trace functions of elements in $H$, but one does not always have a nice set of elements for which the trace functions should be computed.  Moreover, properties of the trace functions that were used in proving results for groups are not true for the elements of the Hopf algebra $H$.  For example, in the group case we showed that if $Tr_A(g,t) = Tr_A(1_G, t)$ then $g = 1_G$ (\cite[Proposition 1.8]{KKZ1}).  However, in the Hopf algebra case, one can add the difference of any two elements with the same trace functions without changing the trace function of an element, so we do not have the strong uniqueness of trace functions that we had for groups.  Characterizing reflection Hopf algebras for AS regular algebras remains an interesting unsolved problem.

\begin{question}
For a Hopf algebra $H$ acting on an AS regular algebra $A$, when is $H$ a reflection Hopf algebra for $A$?
\end{question}

\section{Artin-Schelter Gorenstein subrings of invariants}

Artin-Schelter regular invariant subrings occur under only very special circumstances, but, as H. Bass has noted, Gorenstein rings are ubiquitous. and many of the interesting fixed subrings are Gorenstein rings.  Twenty years after the Shephard-Todd-Chevalley Theorem was proved, Watanabe \cite{W1} showed that if $G$ is a  finite subgroup of ${\rm SL}_n(k)$, then $k[x_1, \dots, x_n]^G$ is a Gorenstein ring, and, in  \cite{W2} he showed that the converse is true if $G$ contains no (classical) reflections.  In our setting, where $A$ is an AS regular algebra, a reasonable generalization of the condition that $A^G$ is a Gorenstein ring is that $A^G$ is an AS Gorenstein algebra.  Next, one must generalize the notion of ``determinant equal to 1".  This generalization was accomplished by P.  J{\o}rgensen and J. Zhang  \cite{JoZ} with their introduction of the notion of the homological determinant of a graded automorphism $g$ of $A$.

The homological determinant $\hdet$ is a group homomorphism $$\hdet: \quad \Aut(A) \rightarrow k^\times$$ that arises in local cohomology; in the case that $A=k[x_1, \dots, x_n]$ it is the determinant (or its inverse, depending upon how $G$ acts on $A$).  The original definition of $\hdet$ is given in \cite[Section 2]{JoZ}, and, in Definition \ref{defhomdet} below, we will give the general definition in the context of Hopf algebra actions.  Fortunately, in many circumstances  $\hdet$ can be computed without using the definition. 
When $A$ is AS regular, the conditions of the following theorem are
satisfied by \cite[Proposition 3.3]{JiZ} and \cite[Proposition 5.5]{JoZ},
and $\hdet g$ can be computed from the trace function of $g$, using the
following result.

\begin{lemma}
\label{xxlem1.4}
\label{hdetbytrace}
\cite[Lemma 2.6]{JoZ}.
Let $A$ be noetherian and AS Gorenstein, and let $g\in \Aut(A)$.
If $g$ is $k$-rational in the sense of \cite[Definition 1.3]{JoZ},
then the rational function ${ Tr}_A(g,t)$ has the form
$${ Tr}_A(g,t) = (-1)^n (\hdet g)^{-1} t^{-\ell}
+ {\rm higher~terms}$$
when it is written as a Laurent series in $t^{-1}$.
\end{lemma}

Our results have shown that the condition that all elements of the group have homological determinant equal to 1 plays the role in our noncommutative setting that the condition that the group is a subgroup of ${\rm SL}_n(k)$ plays in classical invariant theory.
Using homological determinant to replace the usual determinant, J{\o}rgensen and Zhang proved the following generalization of Watanabe's Theorem.
\begin{theorem}\label{watanabe}\cite[Theorem 3.3]{JoZ}.  If $G$ is a finite group of graded automorphisms of an AS regular algebra $A$ with $\hdet(g) = 1$ for all $g \in G$ then $A^G$ is AS Gorenstein.
\end{theorem}

In the classical case, the symmetric group  $S_n$ acting as permutations, is a reflection group for $k[x_1, \dots, x_n]$, but we have already seen (Example \ref{transposition}) that this is not the case for $k_{-1}[x,y]$.  However,  next we note that if $A= k_{-1}[x_1, \dots, x_n]$ (where for each $i \neq j $ we have the relations $x_j x_i = -x_i x_j$), and  if $S_n$ acts on $A$ as permutations, then all subgroups $G$ of $S_n$ have trivial homological determinant, so produce AS Gorenstein invariant subrings. It follows that the fixed subring in Example \ref{transposition} is AS Gorenstein.

\begin{example} \cite[Theorem 5.1]{KKZ5}.
 Let  ${g}$ be a 2-cycle and {$A = k_{-1}[x_1 \ldots, x_n]$}  then
$$Tr_{{A}}({g},t)  = \frac{1}{(1+t^2)(1-t)^{n-2}}$$
$$= (-1)^n \frac{1}{t^n} + \text{ higher terms }$$
so hdet ${g} = 1$,   and hence by Theorem \ref{watanabe}, for ALL groups {$G$} of $n \times n$ permutation matrices, ${A}^{G}$ is
AS Gorenstein.     This, of course, is not true for permutation actions on  a commutative polynomial ring -- e.g.
${k[x_1, x_2, x_3, x_4]}^{\langle (1,2,3,4) \rangle}$
is not Gorenstein,    while
${k_{-1}[x_1, x_2, x_3, x_4]}^{\langle (1,2,3,4) \rangle}$
is AS Gorenstein.
\end{example}
The invariants  of $ k_{-1}[x_1, \dots, x_n]$ under permutation actions are studied in detail in \cite{KKZ5}, producing an interesting contrast to the classical case.  As one example, these groups of permutations contain no reflections of $k_{-1}[x_1, \dots, x_n]$ \cite[Lemma1.7 (4)]{KKZ5}, and so are ``small groups", while the permutation representation of  $S_n$ is a classical reflection group.

A theorem of R. Stanley \cite[Theorem 4.4]{Sta} states that the fixed subring $B= k[x_1, \dots,x_n]^G$ is Gorenstein if and only if the Hilbert series of $B$ satisfies the functional equation
$H_B(t) = \pm t^{-m} H_B(t^{-1})$ for some integer $m$.   J{\o}rgensen and Zhang extended that result to the more general setting of finite groups acting on AS regular algebras.
\begin{theorem} \cite[Proposition 3.8]{KKZ2}. \label{stanley} Let $A$ be an AS regular algebra that satisfies a polynomial identity, and $G$ be a finite group of graded automorphisms of $A$. Then $B= A^G$ is AS Gorenstein if and only if the Hilbert series of $B$ satisfies the functional equation $H_B(t) = \pm t^{-m} H_B(t^{-1})$ for some integer $m$. 
\end{theorem}
Theorem \ref{stanley} is also true under more general (but technical) conditions (see \cite{JoZ} for details).

The homological (co)determinant and Theorems \ref{watanabe} and \ref{stanley} were extended to actions by semisimple Hopf algebras in \cite{KKZ2}.
In \cite{CKWZ1} the homological (co)determinant is defined for any finite dimensional Hopf algebra.  Let $A$ be AS Gorenstein of injective dimension $d$, and let $H$ be a finite dimensional Hopf algebra acting on $A$.  Let $\mathfrak{m}$ denote the maximal graded ideal of $A$ consisting of all elements of positive degree, and let $H_{\mathfrak{m}}^d(A)$ be the $d$-th local cohomology of $A$ with respect to $\mathfrak{m}$.  The lowest degree nonzero homogeneous component of $H_{\mathfrak{m}}^d(A)$ is 1-dimensional; let $\mathfrak{e}$ be a basis element.  Then there is an algebra homomorphism $\eta: H \rightarrow k$ such that the right $H$-action on $H_{\mathfrak{m}}^d(A)^*$ is given by $\eta(h)\mathfrak{e}$ for all $h \in H$. 
\begin{definition} \cite[Definition 3.3]{KKZ2} and \cite[Definition 1.7]{CKWZ1}. \label{defhomdet} Retaining the notation above, 
the composition $\eta \circ S: H \rightarrow k$ is called the {\it homological determinant of the $H$-action on $A$}, and is denoted by $\hdet_H A$.  We say that {\it $\hdet_H A$ is trivial} if $\hdet_H A = \epsilon$, the counit of $H$. 
\end{definition}
Dually, if $K$ coacts on $A$ on the right, then $K$ coacts on $k\mathfrak{e}$ and $\rho(\mathfrak{e}) = \mathfrak{e} \otimes {\rm D}^{-1}$ for some grouplike element ${\rm D}$ in $K$.
\begin{definition}\cite[Definition 6.2]{KKZ2} and \cite[Definition 1.7]{CKWZ1}. Retaining the notation above, 
the {\it homological codeterminant of the $K$-coaction on $A$} is defined to be hcodet$_K A = {\rm D}$.  We say that {\it  hcodet$_K(A)$ is trivial} if hcodet$_K A = 1_K$.
\end{definition}
The homological determinant $\hdet_H A$ is trivial if and only if the homological codeterminant hcodet$_{H^0} A$ is trivial (\cite[Remark 6.3]{KKZ2}).

Watanabe's Theorem was proved for semisimple Hopf actions on AS regular algebras in 
\cite[Theorem 3.6]{KKZ2}; it was also shown for all finite dimensional Hopf actions on AS regular algebras of dimension 2  with trivial homological determinant in \cite[Proposition 0.5]{CKWZ1}.

\begin{theorem}\label{hopfwat} \cite[Theorem 3.6]{KKZ2}.
Let $H$ be a semisimple Hopf algebra acting an AS regular algebra $A$ with trivial homological determinant.  Then $A^H$ is AS Gorenstein.
\end{theorem}
A partial converse to Theorem \ref{hopfwat} is given in \cite[Theorem 4.10]{KKZ2}; in particular if $G$ is a finite group containing no reflections (i.e. a ``small group"), then $A^G$ is AS Gorenstein if and only $\hdet_G(A)$ is trivial, recovering the classical result in \cite{W2}.

If $G$ is a finite subgroup of ${\rm SL}_n(k)$ acting on $A = k[x_1, \dots, x_n]$ then $G$ contains no reflections (since when $g$ is a classical reflection, $\det(g) = $ a root of unity $\neq 1$) so $A^G$ is not a polynomial ring.  We obtain a similar result for the homological determinant.

\begin{theorem} \label{trivialnotregular} \cite[Theorem 2.3]{CKWZ1}. Let $H$ be a semisimple Hopf algebra, and $A$ be a noetherian AS regular algebra equipped with an $H$-algebra action.  If $A^H \neq A$ and the $\hdet_H(A)$ is trivial, then $A^H$ is not AS regular.
\end{theorem}

Since the $\hdet$ is a homomorphism into $k$, it follows from Theorem \ref{trivialnotregular} that if $G$ is a group with $[G,G] = G$ (e.g. if $G$ is a simple group) then $A^G$ will never be AS regular.  Hence such groups are never reflection groups for some AS regular algebra.

We also obtain a version of Stanley's Theorem for semisimple Hopf actions.

\begin{theorem} \cite[Proposition 3.8]{KKZ2}. \label{hopfstanley} Let $A$ be an AS regular algebra that satisfies a polynomial identity, and $H$ be a semisimple Hopf algebra acting on $A$. Then $B= A^G$ is AS Gorenstein if and only if the Hilbert series of $B$ satisfies the functional equation $H_B(t) = \pm t^{-m} H_B(t^{-1})$ for some integer $m$. 
\end{theorem}

In 1884 Felix Klein (\cite{Kl1} \cite{Kl2} \cite{Su}) classified the finite subgroups of ${\rm SL}_2(k)$ and calculated the invariants $k[x,y]^G$, the ``Kleinian singularities", that are important in commutative algebra, algebraic geometry, and representation theory.  These rings of invariants are hypersurfaces in $k[x,y,z]$, and the singularity is of type A, D, or E corresponding to  the type of the McKay quiver of the irreducible representations of the group $G$. The paper \cite{CKWZ1} begins the analogous project for any AS regular algebra of dimension 2, finding all finite dimensional Hopf algebras $H$ that act on $A$, with our standing assumptions (\ref{standing}), and having trivial homological determinant; such a Hopf algebra is called a {\it quantum binary polyhedral group}.  The AS regular algebras of dimension 2 generated in degree one are isomorphic to:
$$k_J[u,v]:=k\langle u,v\rangle/ (vu-uv-u^2)\;\;  {\text{ or}}$$
$$\qquad k_q[u,v]:=k\langle u,v\rangle/(vu-quv).$$
The groups that act on one of these algebras with trivial homological determinant are the cyclic groups, the symmetric group $S_2$, the classical binary polyhedral groups, as well as the dihedral groups (which classically are reflection groups).  
The additional semisimple Hopf algebras that occur are the dual of the group algebra of the dihedral group of order 8, the duals of various finite Hopf quotients of the coordinate Hopf algebra $\mathcal{O}_q({\rm SL}_2(k))$, Hopf algebras that have been studied by \cite{BN}, \cite{Ma2}, \cite{Mu}, and \cite{Ste}.  
In addition, there are actions by non-semisimple Hopf algebras:  the dual of the generalized Taft Hopf algebras $T_{q,\alpha,n}^\circ$, and Hopf algebras whose duals are extensions of the duals of various group algebras by the duals of certain quantum groups.  The table below (reproduced from \cite[Table 1]{CKWZ1}) gives the corresponding AS regular algebras $A$ of dimension 2 and the finite dimensional Hopf algebras $H$ acting on $A$ with trivial homological determinant.\\

\noindent \underline{Notation.} [$\tilde{\Gamma}$, $\Gamma$, $C_n$, $D_{2n}$] 
Let $\tilde{\Gamma}$ denote a finite subgroup of ${\rm SL}_2(k)$, $\Gamma$ 
denote a finite subgroup of ${\rm PSL}_2(k)$, $C_n$ denote a cyclic group of 
order $n$, and $D_{2n}$ denote a dihedral group of order $2n$. Let 
$\text{o}(q)$ denote the order of $q$, for $q \in k^{\times}$ a 
root of unity.  We write $A= k(U)/I$, where $U=ku \oplus kv$, and $I$ is the two-sided ideal generated by the relation.\\
\newpage
\centerline{The quantum binary polyhedral groups $H$ and the AS regular algebras $A$ they act upon:}
\[
\begin{array}{|l|l|}
\hline
\text{AS regular algebra $A$ gldim 2} & \text{finite dimensional Hopf algebra(s) $H$ acting on $A$}\\
\hline
\hline
&  \\
k[u,v] & k\tilde{\Gamma}\\
\hline
&\\
k_{-1}[u,v] & kC_n~\text{for}~n\geq 2; \hspace{.2in} kS_2, \hspace{.2in} kD_{2n};
\\
&  (kD_{2n})^{\circ}; \hspace*{.2in} \mc{D}(\tilde{\Gamma})^{\circ} \text{ for } \tilde{\Gamma} \text{ nonabelian}\\
\hline
&\\
k_q[u,v], ~q \text{ root of 1, } \,q^2 \neq 1 & \\
&\\
{\text{ if $U$ non-simple}} & kC_n \text{~for~} n\geq 3; \hspace{.2in} (T_{q, \alpha, n})^{\circ};\\
\text{ if $U$ simple, $o(q)$  odd} & H \text{ with } 1 \to (k\tilde{\Gamma})^{\circ} \to H^{\circ} \to \mf{u}_{q}(\mf{sl}_2)^{\circ} \to 1;\\
\text{ if $U$ simple, $o(q)$ even, }\,  q^4 \neq 1&  H \text{ with } 1 \to (k\Gamma)^{\circ} \to H^{\circ}
\to \mf{u}_{2,q}(\mf{sl}_2)^{\circ} \to 1;\\
\text{ if $U$ simple}, ~q^4 =1 &
\begin{tabular}{l}
\hspace{-.13in}  $H$ \text{ with } $1 \to (k\Gamma)^{\circ} \to H^{\circ}
\to \mf{u}_{2,q}(\mf{sl}_2)^{\circ} \to 1$   \\
 \hspace{-.13in}  $H$ \text{ with } $1 \to (k\Gamma)^{\circ} \to H^{\circ}
\to  \frac{\mf{u}_{2,q}(\mf{sl}_2)^{\circ}}{(e_{12}-e_{21} e_{11}^2)}\to 1$
\end{tabular}
\\
&\\
\hline
&\\
k_q[u,v], ~q \text{ not root 1} & kC_n, n \geq 2\\
\hline
&\\
k_J[u,v] & kC_2 \\
\hline
\end{array}
\]
~\\
\centerline{Table 1  (\cite[Table 1]{CKWZ1})}
~\\

An interesting next question is to determine when a theorem of Auslander is true in the noncommutative setting.  Recall that a group is called {\it small} if it contains no classical reflections; for example, subgroups of ${\rm SL}_2(k)$ are small.
\begin{theorem}{\rm({\bf Auslander's Theorem})} \cite[Proposition 3.4]{Aus} \rm{and} \cite[Theorem 5.15]{LW}.  Let $G$ be a small finite subgroup of ${\rm GL}_n(k)$ acting linearly on $A=k[x_1, \dots, x_n]$. Then the skew group ring $A\#G$ is naturally isomorphic as a graded algebra to 
the endomorphism ring $\Hom_{A^G}(A,A)$.
\end{theorem}
Some generalizations of Auslander's Theorem were proved by I. Mori and K. Ueyama for groups with trivial homological determinant.  They show \cite[Theorem 3.7]{MU} that if  $G$ is ``ample for $A$" in their sense, then $A\#G$ and $\Hom_{(A^G)^{op}}(A,A)$ are isomorphic as graded algebras, and they give a condition \cite[Corollary 3.11 (3)]{MU} that can be checked for the groups with trivial $\hdet$ acting on AS regular algebras of dimension 2.  They relate Auslander's Theorem to Ueyama's notion of graded isolated singularity \cite{U}.  
In \cite{CKWZ2} Auslander's Theorem  is proved when $A$ has dimension 2 and $H$ is  a semisimple Hopf algebra acting on $A$ under hypotheses (\ref{standing}) and trivial $\hdet_A(H)$.  It is conjectured that  Auslander's Theorem holds for noetherian AS regular algebras $A$ in any dimension.

\begin{conjecture}  If $A$ is an AS regular noetherian algebra and $H$ is a semisimple Hopf algebra acting on $A$ under hypotheses (\ref{standing}) with trivial homological determinant, then
 $A\#G$ is naturally isomorphic to  $\Hom_{(A^G)^{op}}(A,A)$  as graded algebras. 
\end{conjecture}

Auslander's theorem was used to relate  finitely generated projective modules over the skew group ring $k[x,y]\#G$ and maximal Cohen-Macaulay modules over $k[x,y]^G$ when $G$ is a finite  subgroup of ${\rm SL}_2(k)$.  Furthermore, a theorem of Herzog \cite{H} states that the indecomposable maximal Cohen-Macaulay $k[x,y]^G$-modules are precisely the indecomposable $k[x,y]^G$ direct summands of $k[x,y]$.  These and other results of the McKay correspondence are explored in
\cite{CKWZ2}  for $A$ a noetherian AS regular algebra of dimension 2 and $H$ a semisimple Hopf algebra acting on $A$ under hypotheses (\ref{standing}) with trivial homological determinant.  We call a graded $A$-module $M$ an {\it initial} $A$-module if it is generated by $M_0$, and $M_i = 0$ for $i < 0$.  Among the results of \cite{CKWZ2} is the following theorem.
\begin{theorem} \cite{CKWZ2} Let $A$ be a noetherian AS regular algebra of dimension 2 and let $H$ be a semisimple Hopf algebra acting on $A$ under hypotheses (\ref{standing}) with trivial homological determinant. Then there is a bijective correspondence between the isomorphism classes of
\begin{enumerate}
\item indecomposable direct summands of $A$ as right $A^H$-modules
\item indecomposable finitely generated, projective, initial left $\Hom_{(A^H)^{op}}(A,A)$-modules
\item indecomposable finitely generated, projective, initial left $A\#H$-modules
\item simple left $H$-modules
\item indecomposable maximal Cohen-Macaulay $A^H$-modules, up to a degree shift.
\end{enumerate}
\end{theorem}
When $A$ is AS regular of dimension 2 and $H$ is a semisimple Hopf algebra acting on $A$ with trivial homological determinant the invariant subalgebras $A^H$ (called ``Kleinian quantum singularities") are all of the form $C/\Omega C$ for $C$ a noetherian AS regular algebra of dimension 3 and $\Omega$ a normal regular element of $C$, and hence can be regarded as hypersurfaces in an AS regular algebra of dimension 3 (see \cite[Theorem 0.1]{KKZ6} and  \cite{CKWZ2}), and the explicit singularity $\Omega$ is given for each case.  In \cite{CKWZ2} it is shown further that the McKay quiver of $H$ is isomorphic to the Gabriel quiver of the $H$-action on $A$, and the quivers that occur are Euclidean diagrams of types $\widetilde{A}, \widetilde{D}, \widetilde{E}, \widetilde{DL}$, and $\widetilde{L}$.

To conclude this section, we note that it is interesting to compare the roles various groups and Hopf algebras play in the invariant theory of 
$k[x_1, \dots, x_n]$ and  $k_{-1}[x_1, \dots, x_n]$.  In Table 2 we give the classical reflections groups and finite subgroups of ${\rm SL}_n(k)$, and the analogous groups and Hopf algebras for $k_{-1}[x_1, \dots, x_n]$.  Here we use the notation $H_8$ for the Kac-Puljutkin algebra, $A_{4n}$ and $B_{4n}$ as in \cite{Ma2}, $A(\widetilde{\Gamma})$ and $B(\widetilde{\Gamma)}$ as in \cite{BN}.  We notice that some of the same groups play different roles in the two contexts.  For example, the dihedral groups are classical reflection groups, but can act with trivial $\hdet$ on  $k_{-1}[x_1, \dots, x_n]$.  The binary dihedral groups are subgroups of ${\rm SL}_2(k)$ but reflections groups for $k_{-1}[x,y]$.\\
\newpage
~\\
\begin{tabular}{|c|c|c|}
\hline
& $A= \mathbb{C}[x_1, \cdots, x_n]$ & $A=\mathbb{C}_{-1}[x_1, \cdots, x_n]$\\
\hline 
&&\\
Reflection Group for $A$ & &\\ &&\\$n= 2$ &  $D_{2n}=G(n,n,2)$  & $ BD_{4n},$ \\

&& $H_8 =B_8$,\\
&& $A_{4m}$ ($m$ odd),  $ B_{4m}^\circ =B_{4m}$\\
&&\\
$n=3$ & &  $A_{12}, S_4$ (rotations of cube)\\
&&\\
Any $n$ & $C_n$, $S_n$, $G(m,p,n)$ & $C_n$\\
&(34 Exceptional &\\
&for various $n$)&\\
&&\\
\hline
&&\\
Special Linear for $A$ &&\\&&\\
$n=2 $ & $C_n, BD_{4n}, \widetilde{\mathcal{T}}, \widetilde{\mathcal{O}}, \widetilde{\mathcal{I}} $&  $C_n$, $D_{2n}$, $S_2$\\
& &$A_{4m}^\circ, B_{4m}^\circ =B_{4m}, k D_{2n}^\circ$,\\
&& $A(\widetilde{\mathcal{T}})^\circ$, $B(\widetilde{\mathcal{T}})^\circ$, $B(\widetilde{\mathcal{I}})^\circ$,\\
&& $A(\widetilde{\mathcal{O}})^\circ$, $B(\widetilde{\mathcal{O}})^\circ$\\
&&\\
Any $n$ &Finite subgroups of ${\rm SL}_n(k)$& $S_n$ and all subgroups\\
&&\\
\hline
\end{tabular}

\vspace*{.2in}
\centerline{Table 2}
\section{Complete intersection subrings of invariants}
Gorenstein commutative rings can have pathological properties, but
a well-behaved class of Gorenstein commutative rings is the class of graded complete intersections, i.e. the rings  of the form $k[x_1, \dots, x_n]/(f_1, \dots, f_m)$ where $f_1, \dots, f_m$ is a regular sequence of homogeneous elements in $k[x_1, \dots, x_n]$. 
When $A$ is a commutative polynomial ring over $\mathbb{C}$, the 
problem of determining which finite groups $G$ have
the property that $A^G$ is a complete intersection
was solved by N.L. Gordeev \cite{G2} (1986) and (independently) by
H. Nakajima \cite{N2}, \cite{N3} (see the survey \cite{NW}) (1984). 
A key result in this classification  is the theorem of Kac-Watanabe \cite{KW} and Gordeev \cite{G1} that provides a necessary condition: 
if the fixed subring $k[x_1,\cdots,x_n]^G$ (for any finite subgroup
$G\subset GL_n(k)$) is a complete intersection, then $G$ is
generated by bireflections (i.e., elements $g\in GL_n(k)$ such that
$\rank (g-I)\leq 2$ -- i.e. all but two eigenvalues of $g$ are 1). However, the condition that $G$ is generated by bireflections is not sufficient for $k[x_1,\cdots,x_n]^G$ to be a complete intesection.

A first problem in generalizing these results to our setting is that there is not an 
established notion of a complete intersection for noncommutative rings.  
In the commutative graded case a connected graded algebra $A$ is
a complete intersection if one of the following four
equivalent conditions holds \cite[Lemma 1.8]{KKZ4} (which references well-known results from \cite{BH}, \cite{FHT}, \cite{FT} \cite{Gu} \cite{Ta}).\\
\begin{enumerate}
\item[(cci$^\prime$)]
$A\cong k[x_1, \dots, x_d]/(\Omega_1, \cdots, \Omega_n)$, where
$\{\Omega_1, \dots, \Omega_n\}$ is a regular sequence of homogeneous
elements  in $k[x_1, \dots, x_d]$ with $\deg x_i>0$.
\item[(cci)]
$A\cong C/(\Omega_1, \cdots, \Omega_n)$, where $C$ is a noetherian
AS regular algebra and  $\{\Omega_1, \dots, \Omega_n\}$ is a
regular sequence of normalizing homogeneous elements in $C$.
\item[(gci)]
The $\Ext$-algebra $E(A):=\bigoplus_{n=0}^\infty \Ext^n_A(k,k)$ of
$A$ has finite Gelfand-Kirillov dimension.
\item[(nci)]
The $\Ext$-algebra $E(A)$ is noetherian.
\end{enumerate}
~\\

In \cite{KKZ4} we proposed calling a connected graded ring a {\it cci, gci, or nci} if the respective condition above holds for $A$; we called $A$ a {\it hypersurface} if it is a cci when $n=1$ (i.e. of the form $C/(\Omega)$, where $C$ is a noetherian AS regular algebra and $\Omega$ is a regular, normal element of $C$).
In the noncommutative case, unfortunately, the conditions (cci), (gci) and (nci) are
not all equivalent, nor does (gci) or (nci) force $A$ to be
Gorenstein \cite[Example 6.3]{KKZ4}, making it unclear which property 
to use as the proper
generalization of a commutative complete intersection. A direct
generalization to the noncommutative case is condition (cci) which
involves considering regular sequences in {\it any} AS regular
algebra (in the commutative case the only AS regular algebras are
the polynomial algebras), and several researchers have taken an approach to complete intersections that uses regular
sequences. Though the condition (cci) seems to be a good definition
of a noncommutative complete intersection, there are very few tools
available to work with condition (cci), except for explicit
construction and computation, and it is not easy to show condition
(cci) fails, since one needs to consider regular sequences in {\it
any} AS regular algebra.

One relation between these properties that holds in the noncommutative setting is given in the following theorem.
\begin{theorem}\cite[Theorem 1.12(a)]{KKZ4}.
\label{xxthm0.1}  Let $A$ be a
connected graded noncommutative algebra.
  If $A$ satisfies (cci),
then it satisfies (gci).
\end{theorem}
\cite[Example 6.3]{KKZ4} shows that even both (gci) and (nci) together do not 
imply (cci), and \cite[Example 6.2]{KKZ4} shows that (gci) does not imply (nci).

 The  Hilbert series of a commutative complete intersection is a quotient of cyclotomic 
polynomials; we call a noncommutative ring whose Hilbert series has this property {\it cyclotomic}.  A commutative complete intersection is also a Gorenstein ring;  we call $A$ {\it cyclotomic Gorenstein} if it is cyclotomic and AS Gorenstein
\cite[Definition 1.9]{KKZ4}.
\begin{theorem}\cite [Theorem 1.12(b, c)]{KKZ4}. \label{notcyclotomic}
If $A$ satisfies (gci) or (nci), and if the Hilbert series of $A$ is a rational 
function $p(t)/q(t)$ for some coprime integral polynomials $p(t), q(t) 
\in \mathbb{Z}[t]$ with $p(0) = q(0) =1$, then $A$ is cyclotomic.
\end{theorem}
In \cite[Section 2]{KKZ4} we show that certain AS Gorenstein Veronese algebras are not cyclotomic, and hence by Theorem \ref{notcyclotomic} these algebras satisfy none of our conditions for a complete intersection.

In our noncommutative invariant theory context we have produced some examples of invariant subrings that are cci algebras.  Classically the invariants of $k[x_1, \dots, x_n]^{S_n}$ under the permutation representation form a polynomial ring, but, as we noted, this is not the case in $k_{-1}[x_1, \dots, x_n]^{S_n}$.  Under the alternating subgroup $A_n$ of these permutation matrices, the invariants $k[x_1, \dots, x_n]^{A_n}$ are a hypersurface in an $n+1$ dimensional polynomial ring.  In \cite{KKZ5} we show that the two invariant subrings of the skew polynomial rings, $k_{-1}[x_1, \dots, x_n]^{S_n}$ and $k_{-1}[x_1, \dots, x_n]^{A_n}$, are each a cci, and we provide generators for the subring of invariants in each case.
We return to Example \ref{transposition}.

\begin{example} 
\label{transagain}  \cite[Remark 2.6]{KKZ6}. As in Example \ref{transposition} let $S_2$ act on $A = k_{-1}[x,y]$ by interchanging $x$ and $y$.  One set of generators for the fixed subring $A^{S_2}$ is $X:= x+y$ and $Y:=(x-y)(xy)$.  These elements generate a down-up algebra $C$ (an AS regular algebra of dimension 3) with relations
$$YX^2 = X^2Y \;\;\; \text{ and } \;\;\; Y^2X = XY^2,$$
and $A^{S_2} \cong C/\langle \Omega \rangle$, where $\Omega:= Y^2 - \frac{1}{4}X^2 (XY+YX)$, a central regular element of $C$, so that $A^{S_2}$ is a hypersurface in $C$.  
\end{example}

To classify the groups that produce complete intersections, one would like to begin by proving the Kac-Watanabe-Gordeev Theorem: that if $A^G$ is a complete intersection (of some kind) then $G$ must be generated by bireflections.  Toward this end we extend the notion of bireflection, as we extended the notion of a classical reflection, using trace functions.
\begin{definition}\cite[Definition 3.7]{KKZ4}.
 Let $A$ be a noetherian connected graded AS
regular algebra of GK-dimension $n$. We call $g\in \Aut(A)$ a {\it bireflection} if its trace function has the form:
$$Tr_A(g,t) = \frac{1}{(1-t)^{n-2} q(t)}$$
where $q(t)$ is an integral polynomial with $q(1) \neq 0$ (i.e. $1$ is a pole of order $n-2$).  We call it a {\it classical bireflection} if all but two of its eigenvalues are 1.
\end{definition}
The following example suggests that this notion of bireflection based on the  trace function may be useful; in this example the fixed subring is a commutative complete intersection, so it satisfies all the equivalent conditions (cci), (nci), and (gci).
\begin{example} \cite[Example 6.6]{KKZ4}.
${A=k_{-1}[x,y,z]}$ is AS regular of dimension 3,  and the automorphism 
\[{g = \matthree{0 & -1 & 0\\1& 0 &0\\0 & 0 & -1}}\]
acts on it.
The eigenvalues of ${g}$ are $-1, i ,-i$ so ${g}$ is not a classical bireflection.
However, $Tr_{A}({g},t) = 1/((1+t)^2(1-t)) = -1/t^3 +$ higher degree terms and ${g}$  is
a bireflection with hdet ${g} =1$. The fixed subring is
\[{A}^{\langle g \rangle} \cong \frac{k[X,Y,Z,W]}{\langle W^2-(X^2+4Y^2)Z\rangle},\]
a commutative complete intersection.
\end{example}
In the context of permutation actions on $A=k_{-1}[x_1, \dots, x_n]$ we have proved the converse of the
Kac-Watanabe-Gordeev Theorem (a result which is not true in the case $A=k[x_1, \dots, x_n]$).

\begin{theorem}\cite[Theorem 5.4]{KKZ5}. \label{permbi}
If $G$ is a subgroup of $S_n$, represented as permutations of $\{x_1, \dots, x_n\}$, and if $G$ is generated by bireflections (defined in terms of the trace functions), then $k_{-1}[x_1, \dots, x_n]^G$ is a cci.
\end{theorem}
 We conjecture that the Kac-Watanabe-Gordeev Theorem is also true in this context, and we have verified it for $n \leq 4$.

In dimension 2, by  \cite[Theorem 0.1]{KKZ6}
all AS Gorenstein invariant subrings under the actions of finite groups are hypersurfaces in AS regular algebras of dimension 3, and all
automorphisms of finite order are trivially bireflections, and hence 
the first interesting case of the Kac-Watanabe-Gordeev Theorem is in 
dimension 3, so that it is natural to investigate generalizations of 
this theorem for down-up algebras. 

Down-up algebras were defined by Benkart and Roby \cite{BR} 
in 1998 as a tool to study 
the structure of certain posets. 
Noetherian graded down-up algebras $A(\alpha, \beta)$ form a class of AS regular algebras 
of global dimension 3 that are generated in degree 1 by two 
elements $x$ and $y$, with two cubic relations:
$$ y^2x = \alpha yxy + \beta xy^2 \text{ and } yx^2 = \alpha xyx + \beta x^2y$$
for scalars $\alpha, \beta \in k$ with $\beta \neq 0$.
These algebras are not Koszul, but 
they are (3)-Koszul.  Their graded automorphism groups, which depend 
upon the parameters $\alpha$ and $\beta$, were computed in \cite{KK}, 
and are sufficiently rich to provide many non-trivial examples (e.g. 
in two cases the automorphism group is the entire group ${\rm GL}_2(k)$). 
However, it follows from \cite[Proposition 6.4]{KKZ1} that these algebras have no reflections, so all finite
subgroups are ``small", and hence from  \cite[Corollary 4.11]{KKZ2}  $A^G$ is AS Gorenstein if and only if
$\hdet $ is trivial.  Noetherian graded down-up algebras satisfy the following version of the 
Kac-Watanabe-Gordeev Theorem.

\begin{theorem}\cite[Theorem 0.3]{KKZ6}.
 Let $A$ be a graded noetherian down-up
algebra and $G$ be a finite subgroup of $\Aut(A)$. Then the
following are equivalent.
\begin{enumerate}
\item[(C1)]
$A^G$ is a gci.
\item[(C2)]
$A^G$ is cyclotomic Gorenstein and $G$ is generated by
bireflections.
\item[(C3)]
$A^G$ is cyclotomic Gorenstein.
\end{enumerate}
\end{theorem}

In many of the cases for $A$ a noetherian graded down-up algebra the fixed algebras $A^G$ are shown to be a cci, and it is an open question whether that is always the case.

It would be interesting to study the relation of these conditions for other classes of 3-dimensional AS regular algebras, and we have work in progress on actions of groups with trivial homological determinant acting on the generic 3-dimensional Sklyanin algebra.

\section{Related research directions}
In this section we briefly sketch some related directions of research, many of which contain open questions.\\
~\\
\noindent
A. {\bf Degree bounds.}   When computing invariant subrings, it is very useful to have an upper bound on the degrees of the algebra generators of the fixed subring.  In 1916  Emmy Noether \cite{No}
proved that $|G|$, the order of the group $G$, is an upper bound on the degrees of the algebra generators of $k[x_1, \dots, x_n]^G$, for any finite group $G$, when $k$ is a field of  characteristic zero.  The Noether upper bound does not always hold in characteristic p (see e.g. \cite[Example 3.5.7 (a) p. 94]{DK}); the survey paper \cite{Ne} is a good introduction to the problem of finding upper bounds on the degrees of the algebra generators of $k[x_1, \dots, x_n]^G$.  We have seen  (Example \ref{transposition}) that the Noether upper bound does not always hold in our noncommutative setting, for in that example the symmetric group $S_2$ has order 2 and the fixed subring requires a degree 3 generator.  In 2011 
P. Symonds \cite{Sy} proved the upper bound $n(|G|-1)$  (if $n >1$ and $|G| > 1$) on the degrees of the generators of $k[x_1, \dots, x_n]^G$, when $k$ is a field of  characteristic  p; letting $n$ be the number of generators of $A$, this upper bound  also is too small an upper bound for the degrees of the generators in Example \ref{transposition}.
 In the case of permutation
actions on $k[x_1, \dots, x_n]$ there is a smaller upper bound, $\max\{n, {n\choose 2}\}$, on the degrees of the generators of the fixed subring under groups of permutations; this upper bound was proved by M. G{\"o}bel  in 1995 \cite{Go}, and is true in any characteristic.   In \cite[Theorem 2.5]{KKZ5} we prove the bound ${n \choose 2}+\lfloor \frac{n}{2}\rfloor (\lfloor\frac{n}{2}\rfloor+1)$ (which is roughly $3n^2/4$) on the degrees of the generators of $k_{-1}[x_1, \dots, x_n]^G$ for $G$ a group of permutations of $\{x_1, \dots, x_n\}$.  This upper bound follows from a more general upper bound that we state below, a bound  that holds for semisimple Hopf actions on quantum polynomial algebras under certain technical conditions; this upper bound can be viewed as a generalization of Broer's Bound (see \cite[Proposition 3.8.5]{DK})  in the classical case.  In the lemma that follows the field $k$ need not have characteristic zero.
\begin{lemma}[\bf Broer's Bound] \cite[Lemma 2.2]{KKZ5}.
\label{zzlem2.2} Let $A$ be a quantum polynomial algebra of 
dimension $n$ and $C$ an iterated Ore extension
$k[f_1][f_2;\tau_2,\delta_2]\cdots [f_n;\tau_n,\delta_n]$. Assume
that
\begin{enumerate}
\item
$B=A^{H}$ where $H$ is a semisimple Hopf algebra acting on $A$,
\item
$C\subset B\subset A$ and $A_C$ is finitely generated, and
\item
$\deg f_{i}>1$ for at least two distinct $i$'s.
\end{enumerate}
 Then $d_{A^H}$, the maximal degree of the algebra generators of $A^H$, satisifies the inequality"
$$d_{A^H}\leq \ell_C-\ell_A=\sum_{i=1}^n \deg f_i -n,$$
where $\ell_A$ and $\ell_C$ are the
AS indices of $A$ and $C$ respectively.
\end{lemma}
It would be useful to have further upper bounds on the degrees of the generators of the subring of invariants.\\
~\\
B. {\bf Actions on other algebras.} In the work described in this survey thus far  we have assumed that $A$ is a graded algebra, and all actions preserve the grading on $A$.  There is recent work on actions on filtered algebras, such as the Weyl algebras.   Basic properties of this approach were established in \cite{CWWZ}, where it was assumed that the actions preserve the filtration on $A$.

More generally one can consider automorphisms or Hopf actions that may not preserve the filtration on $A$.  Etingof and Walton began a  program to show that in rather  general circumstances a Hopf action must factor through a group action, beginning with \cite[Theorem 1.3]{EW1} that shows that semisimple Hopf actions on  commutative domains must factor through group actions.  In \cite{EW2} actions of finite dimensional Hopf algebras (that are not necessarily semisimple) on commutative domains, particularly when $H$ is pointed of finite Cartan type, are studied;  in this setting there are nontrivial Hopf actions by Taft algebras, Frobenius-Lusztig kernels $u_q({\mathfrak{sl}}_2)$, and Drinfeld twists of some other small quantum groups.  In \cite[Theorem 4.1]{CEW} Cuadra, Etingof and Walton show that if a semisimple Hopf algebra $H$ acts inner faithfully on a Weyl algebra $A_n(k)$ for $k$ an algebraically closed field of characteristic zero, then $H$ is cocommutative;  in this setting they show further \cite[Theorem 4.2]{CEW} that if $H$ is not necessarily semisimple,  but gives rise to a Hopf-Galois extension, then $H$ must be cocommutative.   All of these results have no assumptions regarding preserving a grading or filtration.

Relaxing the noetherian assumption of \cite{CKWZ1}, universal quantum linear group coactions on non-noetherian AS regular algebras of dimension 2 are considered in \cite{WW}.  In another direction, Hopf algebras (including  Taft algebras, doubles of Taft algebras, and $u_q({\mathfrak{sl}}_2)$), that act on certain path algebras,  preserving the path length filtration, are studied in \cite{KiW}.
~\\.  
~\\C. {\bf Nakayama automorphism.} Considering further generalizations of the algebra $A$ on which the group or Hopf algebra acts, let $A$ be a (not necessarily graded) algebra over $k$, and let $A^e = A \otimes A^{op}$ denote the enveloping algebra of $A$.

\begin{definition}\cite[Definition 0.1]{RRZ} and \cite[Section 3.2]{Gi}.

\begin{enumerate}

\item $A$ is called {\it skew Calabi-Yau} (or {\it skew CY} for short) if \begin{enumerate}
\item  $A$ has a projective resolution of finite length in the category
$A^e$-Mod, with every term in the projective resolution finitely generated, and 
\item there is an integer $d$ and an automorphism $\mu$ of $A$ such that 
$$\Ext^i_{A^e}(A,A^e) \cong ~ ^1 A^\mu  \text{ for } i = d, \text{ and } \Ext^i_{A^e}(A,A^e) \cong  0 \text{ if } i \neq d$$
as $A$-bimodules, where $1$ denotes the identity map on $A$.  The map $\mu$ is usually denoted $\mu_A$ and is called the 
{\it Nakayama automorphism of $A$}.
\end{enumerate}
\item \cite[Definition 3.2.3]{Gi} $A$ is called {\it Calabi-Yau} (or {\it CY} for short) if $A$ is skew Calabi-Yau and $\mu_A$ is an inner automorphism of $A$.
\end{enumerate}
\end{definition}
By \cite[Lemma 1.2]{RRZ} 
 if $A$ is a connected graded algebra then $A$ is an AS regular algebra if and only if $A$ is skew CY.

A homological identity is given in \cite[Theorem 0.1]{CWZ} that has been used to show that the Nayakama automorphism plays a role in determining the class of Hopf algebras that can act on a given AS regular algebra.  These techniques were used to show that if a finite dimensional Hopf algebra $H$ acts on $A=k_p[x_1, \dots, x_n]$ (under Hypotheses \ref{standing}) and $p$ is not a root of unity  then $H$ is a group algebra (\cite[Theorem 0.4]{CWZ}); further, if it acts on the 3 or 4-dimensional Sklyanin algebras with trivial $\hdet$ then $H$ is semisimple (\cite[Theorem 0.6]{CWZ}). In \cite{LMZ} the Nakayama automorphism is used to characterize the kinds of Hopf algebras that can act on various families of 3-dimensional AS regular algebras;  in several generic cases it is shown that the Hopf algebra must be a commutative group algebra or (when the $\hdet$ is trivial) the dual of a group algebra.

Further investigation of the Nakayama automorphism is likely to be useful in the study of Hopf actions,  including actions on algebras that are not graded.\\
~\\
D. {\bf  Computing the full automorphism group.} The first step in proving properties of the invariant subring of an algebra $A$ under {\bf any} finite group of automorphisms of $A$ usually is to determine  the complete automorphism group of $A$.  Such computations are notoriously difficult for commutative polynomial rings.  Noncommutative algebras are more rigid than commutative algebras, so sometimes this task is more tractable for noncommutative algebras, even for some PI algebras.  

A sequence of recent papers provide some new techniques for computing the full automorphism group of some algebras, including some filtered algebras whose associated graded algebras are AS regular algebras.
 In \cite{CPWZ1} an invariant, the discriminant of the algebra over a central subring, is defined and used to compute the full automorphism group of some noncommutative algebras, including,  for $n$ even, the filtered algebra (a ``quantum Weyl algebra") $W_n = k \langle x_1, \dots, x_n\rangle$ with relations $x_ix_j + x_j x_i = 1$ for $i > j$,  and its associated graded algebra $k_{-1}[x_1, \dots, x_n]$.  To cite another example, the discriminant is used to show that the full automorphism group  of $B=k_{-1}[x,y]^{S_2}$ (the hypersurface of Example \ref{transposition}) is $k^\times 
\rtimes S_2$ \cite[Example 5.8]{CPWZ1}.  In \cite{CPWZ2}, automorphism groups of tensor products of quantum Weyl algebras and certain skew polynomial rings $k_{q_{i,j}}[x_1, \dots,x_n]$ are computed.  In \cite{CPWZ3} it is shown that, when $n$ is even 
and $n \geq 4$, the fixed subring of $W_n$ under any group of automorphisms of $W_n$ is a filtered AS Gorenstein algebra, but for $n \geq 3$ and odd the full automorphism group of $W_n$ contains a free subalgebra on two (and hence countably many) generators. 
In \cite{CYZ} further results on computing the discriminant are proved, and some applications to Zariski cancellations problems and isomorphism questions of algebras are given (these two areas of application will be explored further in \cite{BZ} and \cite{CPWZ4}).  Further results on the discriminant and its applications remain to be explored.

 This new information about the full automorphism group of many families of noncommutative algebras suggests  many interesting open questions about the structure of the invariant subrings.  The questions considered in this survey can be investigated for actions of ANY finite group of (not necessarily graded) automorphisms of $A$, for larger classes of algebras than AS regular algebras.

\subsection*{Acknowledgments}
The author wishes to thank Chelsea Walton and James Zhang, as well as the referee, for making helpful suggestions on this paper.
Ellen Kirkman was partially supported by the 
Simons Foundation grant no. 208314

\end{document}